\documentclass{amsart}
\usepackage{graphicx}
\usepackage[english]{babel}
\usepackage{amssymb}
\usepackage{graphicx}
\usepackage{natbib}
\usepackage{amsfonts}
\usepackage{amscd}
\usepackage{amsmath}
\usepackage{amsthm}
\usepackage{latexsym}
\usepackage[latin1]{inputenc}
\usepackage[T1]{fontenc}
\usepackage[pdftex,colorlinks,linkcolor=red,citecolor=blue]{hyperref}
\setcitestyle{notesep={; },authoryear,round,coma,aysep={,},yysep={;}}
\newtheorem{Thm}{Theorem}[section]
\newtheorem{Coro}[Thm]{Corollary}
\newtheorem{Lem}[Thm]{Lemma}

\theoremstyle{definition}
\newtheorem{Def}[Thm]{Definition}
\theoremstyle{remark}
\newtheorem{Rem}[Thm]{Remark}
\newtheorem*{pf}{Proof}
\numberwithin{equation}{section}

\usepackage{tikz,xcolor,hyperref}
\definecolor{lime}{HTML}{A6CE39}
\DeclareRobustCommand{\orcidicon}{
	\begin{tikzpicture}
	\draw[lime, fill=lime] (0,0)
	circle [radius=0.16]
	node[white] {{\fontfamily{qag}\selectfont \tiny ID}};
	\draw[white, fill=white] (-0.0625,0.095)
	circle [radius=0.007];
	\end{tikzpicture}
	\hspace{-2mm}
}
\foreach \x in {A, ..., Z}{\expandafter\xdef\csname orcid\x\endcsname{\noexpand\href{https://orcid.org/\csname orcidauthor\x\endcsname}
			{\noexpand\orcidicon}}
}
\begin{document}
\bibliographystyle{abbrvnat}

\title[LDP for a $G$-BSDE with subdifferential operator]{Large deviation principle for a backward
stochastic differential equation driven by $G$-Brownian motion with subdifferential operator}%
\author[A.~S.~Hima]{Abdoulaye Soumana Hima}%
\address{Département de Mathématiques, Université Dan Dicko Dankoulodo de Maradi, BP 465, Maradi, Niger}%
\email[A.~S.~Hima]{abdoulaye.hima@uddm.edu.ne}

\author[I.~Dakaou]{Ibrahim Dakaou\orcidA{}}%
\email[I.~Dakaou]{ibrahim.dakaou@uddm.edu.ne}

\subjclass[2010]{Primary 60F10; Secondary 60H10, 60H30}%
\keywords{Large deviations, contraction principle, backward stochastic differential equation, $G$-Brownian motion, subdifferential operator, variational inequality}%

\begin{abstract}
In this paper, we study a large deviation principle for the solution of a backward stochastic differential equation driven by $G$-Brownian motion with subdifferential operator.
\end{abstract}
\maketitle
\section{Introduction}
The large deviation principle (LDP in short) characterizes the limiting behavior, as $\varepsilon\rightarrow0$, of family of probability measures $\{\mu_{\varepsilon}\}_{\varepsilon>0}$ in terms of a rate function. Several authors have considered large deviations and obtained different types of applications mainly to mathematical physics. General references on large deviations are: \citet{Varadhan1984,Deuschel1989,Dembo1998}.

Let $X^{s,x,\varepsilon}$ be the diffusion process that is the unique solution of the following stochastic differential equation (SDE in short)
\begin{equation}\label{eq0}
X_{t}^{s,x,\varepsilon}=x+\int_{s}^{t}\beta(X_{r}^{s,x,\varepsilon})dr+\sqrt
\varepsilon \int_{s}^{t}\sigma(X_{r}^{s,x,\varepsilon})dW_{r},\; 0\leq s\leq t\leq T
\end{equation}
where $\beta$ is a Lipschitz function defined on $\mathbb{R}^{n}$ with values in $\mathbb{R}^{n}$, $\sigma$ is a Lipschitz function defined on $\mathbb{R}^{n}$ with values in $\mathbb{R}^{n\times d}$, and $W$ is a standard Brownian motion in $\mathbb{R}^{d}$ defined on a complete probability space $(\Omega,\mathcal{F},\mathbb{P})$. The existence and uniqueness of
the strong solution $X^{s,x,\varepsilon}$ of \eqref{eq0} is standard. Thanks to the work of \citet{Freidlin1984}, the sequence $(X^{s,x,\varepsilon})_{\varepsilon > 0}$ converges in probability, as $\varepsilon$ goes to 0, to $(\varphi_{t}^{s,x})_{s\leq t\leq T}$ solution  of the following deterministic equation
\begin{equation*}
\varphi_{t}^{s,x}=x+\int_{s}^{t}\beta(\varphi_{r}^{s,x})dr,\; 0\leq s\leq t\leq T
\end{equation*}
and satisfies a LDP.

\citet{Rainero2006} extended this result to the case of backward stochastic differential equations (BSDEs in short) and \citet{Essaky2008} to BSDEs with subdifferential operator.

\citet{Gao2010} extended the work of \citet{Freidlin1984} to stochastic differential equations driven by $G$-Brownian motion ($G$-SDEs in short). The authors considered the following $G$-SDE: for every $0\leq t\leq T$,
\begin{equation*}
    X_{t}^{x,\varepsilon} =x+\int_{0}^{t}b^{\varepsilon}(X_{r}^{x,\varepsilon})dr+\varepsilon\int_{0}^{t}h^{\varepsilon}(X_{r}^{x,\varepsilon})d\langle B, B\rangle_{r/\varepsilon}+\varepsilon\int_{0}^{t}\sigma^{\varepsilon}(X_{r}^{x,\varepsilon})dB_{r/\varepsilon}
\end{equation*}
and use discrete time approximation to establish LDP for $G$-SDEs.

\citet{Hu2014a} proved the existence and uniqueness of the solutions for BSDEs driven by G-Brownian motion. Moreover, \citet{Hu2014b} showed the comparison theorem, \emph{Feynman-Kac} formula, and \emph{Girsanov} transformation for G-BSDEs and established the probabilistic interpretation for the viscosity solutions of a class of fully nonlinear partial differential equations (PDEs in short).

\citet{Yang2017} proved the existence and uniqueness of a solution for a class of BSDEs driven by $G$-Brownian motion with subdifferential operator  ($G$-MBSDEs in short) and established a probabilistic interpretation for the viscosity solutions of a class of nonlinear variational inequalities.

Recently, \citet{Dakaou2021} established a LDP for $G$-BSDEs. More precisely, the authors considered the following forward-backward stochastic differential equation driven by $G$-Brownian motion: for every $s\leq t\leq T$,
\begin{equation*}
\begin{cases}
&X_{t}^{s,x,\varepsilon} =x+\int_{s}^{t}b(X_{r}^{s,x,\varepsilon})dr+\varepsilon\int_{s}^{t}h(X_{r}^{s,x,\varepsilon})d\langle B, B\rangle_{r}+\varepsilon\int_{s}^{t}\sigma(X_{r}^{s,x,\varepsilon})dB_{r} \\
&Y_{t}^{s,x,\varepsilon}=\Phi(X_{T}^{s,x,\varepsilon})+\int_{t}^{T}f(r,X_{r}^{s,x,\varepsilon},Y_{r}^{s,x,\varepsilon},Z_{r}^{s,x,\varepsilon})dr-\int_{t}^{T}Z_{r}^{s,x,\varepsilon}dB_{r}
\\ &\qquad\qquad+\int_{t}^{T}g(r,X_{r}^{s,x,\varepsilon},Y_{r}^{s,x,\varepsilon},Z_{r}^{s,x,\varepsilon})d\langle B, B\rangle_{r}-(K_{T}^{s,x,\varepsilon}-K_{t}^{s,x,\varepsilon})
\end{cases}
\end{equation*}
They studied the asymptotic behavior of the solution of the backward equation and established a LDP for the corresponding process.

Motivated by the aforementioned works, we aim to establish LDP for $G$-BSDEs with subdifferential operator. More precisely, we consider the following forward-backward stochastic differential equation driven by $G$-Brownian motion with subdifferential operator: for every $s\leq t\leq T$,
\begin{equation*}
\begin{cases}
&X_{t}^{s,x,\varepsilon} =x+\int_{s}^{t}b(X_{r}^{s,x,\varepsilon})dr+\varepsilon\int_{s}^{t}h(X_{r}^{s,x,\varepsilon})d\langle B, B\rangle_{r}+\varepsilon\int_{s}^{t}\sigma(X_{r}^{s,x,\varepsilon})dB_{r} \\
&-dY_{t}^{s,x,\varepsilon}+\partial\Pi(Y_{t}^{s,x,\varepsilon})dt\ni f(t,X_{t}^{s,x,\varepsilon},Y_{t}^{s,x,\varepsilon},Z_{t}^{s,x,\varepsilon})dt-Z_{t}^{s,x,\varepsilon}dB_{t} \\
&\qquad\qquad\qquad\qquad\qquad\qquad+g(t,X_{t}^{s,x,\varepsilon},Y_{t}^{s,x,\varepsilon},Z_{t}^{s,x,\varepsilon})d\langle B\rangle_{t}-dK_{t}^{s,x,\varepsilon} \\
&Y_{T}^{s,x,\varepsilon}=\Phi(X_{T}^{s,x,\varepsilon})
\end{cases}
\end{equation*}
where $K^{s,x,\varepsilon}$ is a decreasing $G$-martingale; $\partial\Pi$ is the subdifferential operator associated with $\Pi$ which is a lower semicontinuous function. The objective of this work is to study the asymptotic behavior of the family $(Y^{s,x,\varepsilon})_{\varepsilon > 0}$ as $\varepsilon$ goes to $0$ and prove that $(Y^{s,x,\varepsilon})_{\varepsilon > 0}$ satisfies a LDP.

The remaining part of the paper is organized as follows. In Section~\ref{Sec2}, we present some preliminaries that are
useful in this paper. Section~\ref{Sec3} is devoted to the large deviations for $G$-SDEs obtained by \citet{Gao2010}.
The large deviations for $G$-MBSDEs are given in Section~\ref{Sec4}.
\section{Preliminaries}\label{Sec2}
We review some basic notions and results about $G$-expectation, $G$-Brownian motion and $G$-stochastic integrals \citep[see][for more details]{Peng2019}.

Let $\Omega$ be a complete separable metric space, and let $\mathcal{H}$ be a linear space of real-valued functions defined on $\Omega$ satisfying: if $X_i\in\mathcal{H}$, $i=1, \ldots, n$, then
\begin{equation*}
\varphi(X_1, \ldots, X_n)\in\mathcal{H}, \quad \forall\varphi\in\mathcal{C}_{l, Lip}(\mathbb{R}^n),
\end{equation*}
where $\mathcal{C}_{l,Lip}(\mathbb{R}^{n})$ is the space of real continuous functions defined on $\mathbb{R}^{n}$ such that for some $C>0$ and $k\in\mathbb{N}$ depending on $\varphi$,
\begin{equation*}
    \vert\varphi (x)-\varphi (y)\vert\leq C(1+\vert x\vert^{k}+\vert y\vert^{k})\vert x-y\vert,\quad \forall x,y\in\mathbb{R}^{n}.
\end{equation*}
\begin{Def}(\emph{Sublinear expectation space}).
A sublinear expectation $\widehat{\mathbb{E}}$ on $\mathcal{H}$ is a functional $\widehat{\mathbb{E}}:\mathcal{H}\longrightarrow \mathbb{R}$ satisfying the following properties: for all $X, Y\in\mathcal{H}$, we have
\begin{enumerate}
  \item Monotonicity: if $X\geq Y$, then $\widehat{\mathbb{E}}[X]\geq \widehat{\mathbb{E}}[Y]$;
  \item Constant preservation: $\widehat{\mathbb{E}}[c]=c$;
  \item Sub-additivity: $\widehat{\mathbb{E}}[X+Y]\leq \widehat{\mathbb{E}}[X]+\widehat{\mathbb{E}}[Y]$;
  \item Positive homogeneity: $\widehat{\mathbb{E}}[\lambda X]=\lambda\widehat{\mathbb{E}}[X]$, for all $\lambda\geq 0$.
\end{enumerate}
$(\Omega, \mathcal{H}, \widehat{\mathbb{E}})$ is called a \emph{sublinear expectation space}.
\end{Def}
\begin{Def}(\emph{Independence}).
Fix the sublinear expectation space $(\Omega ,\mathcal{H},\widehat{\mathbb{E}})$. A random variable $Y\in \mathcal{H}$ is said to be independent of $(X_{1}, X_{2}, \ldots, X_{n})$, $X_{i}\in\mathcal{H}$, if for all $\varphi\in\mathcal{C}_{l, Lip}(\mathbb{R}^{n+1})$,
\begin{equation*}
    \widehat{\mathbb{E}}\left[\varphi (X_{1}, X_{2}, \ldots, X_{n}, Y)\right]=
    \widehat{\mathbb{E}}\left[\widehat{\mathbb{E}}\left[\varphi(x_{1}, x_{2}, \ldots, x_{n}, Y)\right]\big|_{(x_{1}, x_{2}, \ldots, x_{n})=(X_{1}, X_{2}, \ldots, X_{n})}\right].
\end{equation*}
\end{Def}
Now we introduce the definition of $G$-normal distribution.
\begin{Def}(\emph{$G$-normal distribution}).
A random variable $X\in\mathcal{H}$ is called $G$-normally distributed, noted by $X\sim\mathcal{N}(0, [\underline{\sigma}^{2}, \overline{\sigma}^{2}])$, $0\leq\underline{\sigma}^{2}\leq\overline{\sigma}^{2}$, if for any function $\varphi\in\mathcal{C}_{l, Lip}(\mathbb{R})$, the fonction $u$ defined by $u(t,x):=\widehat{\mathbb{E}}[\varphi(x+\sqrt{t}X)],\;\left(t, x\right)\in\left[0,\infty\right)\times\mathbb{R}$, is a viscosity solution of the following $G$-heat equation:
\begin{equation*}\label{heat}
    \partial_{t}u-G\left(D_{x}^{2}u\right)=0,\;u(0,x)=\varphi(x),
\end{equation*}
where
\begin{equation*}
    G(a):=\frac{1}{2}(\overline{\sigma}^{2}a^{+}-\underline{\sigma}^{2}a^{-}).
\end{equation*}
In multi-dimensional case, the function $G(\cdot)$: $\mathbb{S}_{d}\longrightarrow\mathbb{R}$  is defined by
\begin{equation*}
    G(A)=\frac{1}{2}\sup_{\gamma\in\Gamma}\,\textrm{tr}(\gamma\gamma^{\tau}A),
\end{equation*}
where $\mathbb{S}_{d}$ denotes the space of $d\times d$ symmetric matrices and $\Gamma$ is a given nonempty, bounded and closed subset of $\mathbb{R}^{d\times d}$ which is the space of all $d\times d$ matrices.
\end{Def}
Throughout this paper, we consider only the non-degenerate case, i.e., $\underline{\sigma}^{2}>0$.

Let $\Omega:=\mathcal{C}\left([0, \infty), \mathbb{R}\right)$, which equipped with the raw filtration $\mathcal{F}$ generated by the canonical process $(B_t)_{t\geq 0}$, i.e., $B_t(\omega)=\omega_t$, for $(t, \omega)\in[0, \infty)\times\Omega$. Let $\Omega_{T}:=\mathcal{C}\left([0, T], \mathbb{R}\right)$ and let us consider the function spaces defined by
\begin{align*}
  Lip(\Omega_{T})&:=\Big\{\varphi(B_{t_{1}}, B_{t_{2}}-B_{t_{1}}, \ldots, B_{t_{n}}-B_{t_{n-1}}):  n\geq 1,\\
  &\quad\quad\quad\quad 0\leq t_{1}\leq t_{2}\leq\ldots\leq t_{n}\leq T, \varphi\in\mathcal{C}_{l,Lip}(\mathbb{R}^{n})\Big\},\quad \text{for}\quad T>0, \\
  Lip(\Omega)&:=\bigcup_{n=1}^{\infty}Lip(\Omega_{n}).
\end{align*}
\begin{Def}(\emph{$G$-Brownian motion and $G$-expectation}).
On the sublinear expectation space $\left(\Omega, Lip(\Omega), \widehat{\mathbb{E}}\right)$, the canonical process $(B_t)_{t\geq 0}$ is called a $G$-Brownian motion if the following properties are verified:
\begin{enumerate}
  \item $B_{0}=0$
  \item For each $t, s\geq 0$, the increment $B_{t+s}-B_t\sim\mathcal{N}(0, [s\underline{\sigma}^{2}, s\overline{\sigma}^{2}])$ and is independent from $(B_{t_{1}}, \ldots, B_{t_{n}})$, for $0\leq t_1\leq \ldots\leq t_n\leq t$.
\end{enumerate}
Moreover, the sublinear expectation $\widehat{\mathbb{E}}$ is called \emph{$G$-expectation}.
\end{Def}
\begin{Rem}\label{l6}
For each $\lambda>0$, $\left(\sqrt{\lambda}B_{t/\lambda}\right)_{t\geq 0}$ is also a $G$-Brownian motion. This is the \emph{scaling property} of $G$-Brownian motion, which is the same as that of the classical Brownian motion.
\end{Rem}
\begin{Def}(\emph{Conditional $G$-expectation}).
For the random variable $\xi\in Lip(\Omega_T)$ of the following form:
\begin{equation*}
\varphi(B_{t_{1}}, B_{t_{2}}-B_{t_{1}}, \ldots, B_{t_{n}}-B_{t_{n-1}}), \quad \varphi\in\mathcal{C}_{l, Lip}(\mathbb{R}^n),
\end{equation*}
the conditional $G$-expectation $\widehat{\mathbb{E}}_{t_{i}}[\cdot]$, $i=1, \ldots, n$, is defined as follows
\begin{equation*}
    \widehat{\mathbb{E}}_{t_{i}}[\varphi(B_{t_{1}}, B_{t_{2}}-B_{t_{1}}, \ldots, B_{t_{n}}-B_{t_{n-1}})]=\widetilde{\varphi}(B_{t_{1}},
    B_{t_{2}}-B_{t_{1}}, \ldots, B_{t_{i}}-B_{t_{i-1}}),
\end{equation*}
where
\begin{equation*}
    \widetilde{\varphi}\left(x_{1}, \ldots, x_{i}\right)=\widehat{\mathbb{E}}\left[\varphi\left(x_{1}, \ldots, x_{i}, B_{t_{i+1}}-B_{t_{i}}, \ldots, B_{t_{n}}-B_{t_{n-1}}\right)\right].
\end{equation*}
If $t\in(t_{i}, t_{i+1})$, then the conditional $G$-expectation $\widehat{\mathbb{E}}_{t}[\xi]$ could be defined by reformulating $\xi$ as
\begin{equation*}
    \xi=\widehat{\varphi}(B_{t_{1}}, B_{t_{2}}-B_{t_{1}}, \ldots, B_{t}-B_{t_{i}}, B_{t_{i+1}}-B_{t}, \ldots, B_{t_{n}}-B_{t_{n-1}}),
    \quad\widehat{\varphi}\in\mathcal{C}_{l,Lip}(\mathbb{R}^{n+1}).
\end{equation*}
\end{Def}
For $\xi\in Lip(\Omega_T)$ and $p\geq 1$, we consider the norm $\Vert\xi\Vert_{L^p_G}:=\left(\widehat{\mathbb{E}}\Big[\vert\xi\vert^p\Big]\right)^{1/p}$. Denote by $L^p_G(\Omega_T)$ the Banach completion of $Lip(\Omega_T)$ under $\Vert\cdot\Vert_{L^p_G}$. It is easy to check that the conditional $G$-expectation $\widehat{\mathbb{E}}_{t}[\cdot]: Lip(\Omega_T)\longrightarrow Lip(\Omega_t)$ is a continuous mapping and thus can be extended to $\widehat{\mathbb{E}}_{t}[\cdot]: L^p_G(\Omega_T)\longrightarrow L^p_G(\Omega_t)$.
\begin{Def}(\emph{$G$-martingale}).
A process $M=\left(M_{t}\right)_{t\in\left[0, T\right]}$ with $M_{t}\in L_{G}^{1}(\Omega_{t})$, $0\leq t\leq T$, is called a $G$-martingale if for all $0\leq s\leq t\leq T$, we have
\begin{equation*}
    \widehat{\mathbb{E}}_{s}[M_{t}]=M_{s}.
\end{equation*}
The process $M=\left(M_{t}\right)_{t\in\left[0, T\right]}$ is called symmetric $G$-martingale if $-M$ is also a $G$-martingale.
\end{Def}
\begin{Thm}\citep[Representation theorem of $G$-expectation, see][]{Hu2009,Denis2011}.
There exists a weakly compact set $\mathcal{P}\subset\mathcal{M}_1(\Omega_T)$, the set of probability measures on $(\Omega_T, \mathcal{B}(\Omega_T))$, such that
\begin{equation*}
    \widehat{\mathbb{E}}[\xi]:=\sup_{P\in \mathcal{P}}E_{P}[\xi] \quad \text{for all} \quad \xi\in L^1_G(\Omega_T).
\end{equation*}
$\mathcal{P}$ is called a set that represents $\widehat{\mathbb{E}}$.
\end{Thm}
Let $\mathcal{P}$ be a weakly compact set that represents $\widehat{\mathbb{E}}$. For this $\mathcal{P}$, we define the \emph{capacity} of a measurable set $A$ by
\begin{equation*}
    \widehat{C}(A):=\sup_{P\in\mathcal{P}}P(A), \quad A\in\mathcal{B}(\Omega_{T}).
\end{equation*}
A set $A\in\mathcal{B}(\Omega_{T})$ is a polar if $\widehat{C}(A)=0$. A property holds \emph{quasi-surely} (q.s.) if it is true outside a polar set.

An important feature of the $G$-expectation framework is that the quadratic variation $\left\langle B\right\rangle$ of the $G$-Brownian motion is no longer a deterministic process, which is given by
\begin{equation*}
    \left\langle B\right\rangle_{t}:=\lim_{\delta\left(\pi_{t}^{N}\right)
    \rightarrow 0}\sum_{j=0}^{N-1}(B_{t_{j+1}^{N}}-B_{t_{j}^{N}})^{2},
\end{equation*}
where $\pi_{t}^{N}=\{t_{0}, t_{1}, \ldots, t_{N}\}$, $N=1, 2, \ldots$, are refining partitions of $[0, t]$. For all $s, t\geq 0$, $\langle B\rangle_{t+s}-\langle B\rangle_{t}\in[s\underline{\sigma}^2, s\overline{\sigma}^2]$, $q.s.$ \citep[see][]{Peng2019}.

Let $M_{G}^{0}\left(0, T\right)$ be the collection of processes in the following form: for a given partition $\pi_{T}^{N}:=\{t_{0}, t_{1}, \ldots, t_{N}\}$ of $[0, T]$,
\begin{equation}\label{simp}
\eta_{t}\left(\omega\right)=\sum_{j=0}^{N-1}\xi_{j}\left(\omega\right)
\mathbf{1}_{\left[t_{j}, t_{j+1}\right)}(t),
\end{equation}
where $\xi_{i}\in Lip(\Omega_{t_{i}})$, for all $i=0, 1, \ldots, N-1$. For $p\geq 1$ and $\eta\in M_{G}^{0}\left(0, T\right)$, let $\left\Vert\eta\right\Vert_{H_{G}^{p}}:=\left(\widehat{\mathbb{E}}\left[\left(\int_{0}^{T}|\eta_{s}|^{2}ds\right)^{p/2}\right]    \right)^{1/p}$, $\Vert\eta\Vert_{M_{G}^{p}}:=\left(\widehat{\mathbb{E}}\left[\int_{0}^{T}|\eta_{s}|^{p}ds\right]\right)^{1/p}$ and denote by $H_{G}^{p}\left(0,T\right)$, $M_{G}^{p}(0,T)$ the completions of $M_{G}^{0}\left(0,T\right)$ under the norms $\Vert\cdot\Vert_{H_{G}^{p}}$, $\Vert\cdot\Vert_{M_{G}^{p}}$ respectively.

Let $\mathcal{S}_{G}^{0}\left(0,T\right):=\{h(t, B_{t_{1}\wedge t},B_{t_{2}\wedge t}-B_{t_{1}\wedge t}, \ldots, B_{t_{n}\wedge
t}-B_{t_{n-1}\wedge t}): 0\leq t_{1}\leq t_{2}\leq \ldots\leq t_{n}\leq T, ~h\in\mathcal{C}_{b,Lip}(\mathbb{R}^{n+1})\}$, where $\mathcal{C}_{b,Lip}(\mathbb{R}^{n+1})$ is the collection of all bounded and Lipschitz functions on $\mathbb{R}^{n+1}$. For $p\geq 1$ and $\eta\in\mathcal{S}_{G}^{0}\left(0, T\right)$, we set $\left\Vert\eta\right\Vert_{\mathcal{S}_{G}^{p}}:=\left(\widehat{\mathbb{E}}\Big[\sup\limits_{t\in[0, T]}\vert\eta_{t}\vert^{p}\Big]\right)^{1/p}$. We denote by $\mathcal{S}_{G}^{p}\left(0, T\right)$ the completion of $\mathcal{S}_{G}^{0}\left(0, T\right)$ under the norm $\left\Vert\cdot\right\Vert_{\mathcal{S}_{G}^{p}}$.
\begin{Def}
For $\eta\in M_{G}^{0}\left(0, T\right)$ of the form \eqref{simp}, the It\^{o} integral with respect to $G$-Brownian motion is defined by the linear mapping $\mathcal{I}: M_{G}^{0}(0, T)\longrightarrow L_{G}^{2}(\Omega_{T})$,
\begin{equation*}
    \mathcal{I}(\eta):=\int_{0}^{T}\eta_{t}dB_{t}=\sum_{k=0}^{N-1}\xi_{k}(B_{t_{k+1}}-B_{t_{k}}),
\end{equation*}
which can be continuously extended to $\mathcal{I}: H_{G}^{1}(0,T)\longrightarrow L_{G}^{2}(\Omega_{T})$. On the other hand, the stochastic integral with respect to $(\langle B\rangle_t)_{t\geq 0}$ is defined by the linear mapping $\mathcal{Q}: M_{G}^{0}(0, T)\longrightarrow L_{G}^{1}(\Omega_{T})$,
\begin{equation*}
    \mathcal{Q}(\eta):=\int_{0}^{T}\eta_{t}d\langle B\rangle_{t}=\sum_{k=0}^{N-1}\xi_{k}(\langle B\rangle_{t_{k+1}}-\langle B\rangle_{t_{k}}),
\end{equation*}
which can be continuously extended to $\mathcal{Q}: H_{G}^{1}(0,T)\longrightarrow L_{G}^{1}(\Omega_{T})$.
\end{Def}
\begin{Lem}\label{BDG}\citep[BDG type inequality, see][Theorem~2.1]{Gao2009}.
Let $p\geq 2$, $\eta\in H_{G}^{p}(0, T)$ and $0\leq s\leq t\leq T$. Then,
\begin{eqnarray*}
  &&c_{p}\underline{\sigma}^{p}\widehat{\mathbb{E}}
    \left[\left(\int_{0}^{T}|\eta_{s}|^{2}ds\right)^{p/2}\right]  \\
  && \leq
    \widehat{\mathbb{E}}\left[\sup_{0\leq t\leq T}\left\vert\int_{0}^{t}
    \eta_{r}dB_{r}\right\vert^{p}\right]\leq
    C_{p}\overline{\sigma}^{p}\widehat{\mathbb{E}}
    \left[\left(\int_{0}^{T}|\eta_{s}|^{2}ds\right)^{p/2}\right],
\end{eqnarray*}
where $0<c_{p}<C_{p}<\infty $ are constants independent of $\eta$, $\underline{\sigma}$ and $\overline{\sigma}$.
\end{Lem}
For $\xi\in Lip(\Omega_{T})$, let
\begin{equation*}
    \mathcal{E}(\xi):=\widehat{\mathbb{E}}\left(\sup_{t\in[0, T]}\widehat{\mathbb{E}}_{t}[\xi]\right).
\end{equation*}
$\mathcal{E}$ is called the $G$-evaluation. For $p\geq 1$ and $\xi\in Lip(\Omega_{T})$, define
\begin{equation*}
    \Vert\xi\Vert_{p, \mathcal{E}}:=\left(\mathcal{E}[\vert\xi\vert^{p}]\right)^{1/p}
\end{equation*}
and denote by $L_{\mathcal{E}}^{p}(\Omega_{T})$ the completion of $Lip(\Omega_{T})$ under the norm $\Vert\cdot\Vert_{p, \mathcal{E}}$.
\begin{Thm}\label{thm2}\citep[See][]{Song2011}.
For any $\alpha\geq 1$ and $\delta>0$, we have $L_{G}^{\alpha+\delta}(\Omega_{T})\subset L_{\mathcal{E}}^{\alpha}(\Omega_{T})$.
More precisely, for any $1<\gamma<\beta:=(\alpha+\delta)/\alpha$, $\gamma\leq 2$ and for all $\xi\in Lip(\Omega_{T})$, we have
\begin{equation*}\label{}
    \widehat{\mathbb{E}}\Big[\sup_{t\in[0, T]}\widehat{\mathbb{E}}_{t}[\vert\xi\vert^{\alpha}]\Big]\leq
    C\Big\{(\widehat{\mathbb{E}}[\vert\xi\vert^{\alpha+\delta}])^{\alpha/(\alpha+\delta)}+
    (\widehat{\mathbb{E}}[\vert\xi\vert^{\alpha+\delta}])^{1/\gamma}\Big\},
\end{equation*}
where
\begin{equation*}
    C=\frac{\gamma}{\gamma-1}(1+14\sum_{i=1}^{\infty}i^{-\beta/\gamma}).
\end{equation*}
\end{Thm}
\begin{Rem}\label{rm1}
By $\frac{\alpha}{\alpha+\delta}<\frac{1}{\gamma}<1$, we have
\begin{equation*}\label{}
    \widehat{\mathbb{E}}\Big[\sup_{t\in[0, T]}\widehat{\mathbb{E}}_{t}[\vert\xi\vert^{\alpha}]\Big]\leq
    2C\Big\{(\widehat{\mathbb{E}}[\vert\xi\vert^{\alpha+\delta}])^{\alpha/(\alpha+\delta)}+
    \widehat{\mathbb{E}}[\vert\xi\vert^{\alpha+\delta}]\Big\}.
\end{equation*}
Set
\begin{equation*}
    C_{1}=2\inf\Big\{\frac{\gamma}{\gamma-1}(1+14\sum_{i=1}^{\infty}i^{-\beta/\gamma}):1<\gamma<\beta, \gamma\leq 2\Big\},
\end{equation*}
then
\begin{equation}\label{eq8}
    \widehat{\mathbb{E}}\Big[\sup_{t\in[0, T]}\widehat{\mathbb{E}}_{t}[\vert\xi\vert^{\alpha}]\Big]\leq
    C_{1}\Big\{(\widehat{\mathbb{E}}[\vert\xi\vert^{\alpha+\delta}])^{\alpha/(\alpha+\delta)}+
    \widehat{\mathbb{E}}[\vert\xi\vert^{\alpha+\delta}]\Big\},
\end{equation}
where $C_{1}$ is a constant only depending on $\alpha$ and $\delta$.
\end{Rem}
\begin{Lem}\label{lemx1}\citep[See][]{Hu2014a}. Let $X\in\mathcal{S}_{G}^{\alpha }(0,T)$ for some $\alpha >1$ and $\alpha ^{\ast }=\dfrac{\alpha }{\alpha -1}$. Assume that $K^{j},~\ j=1,2$, are two decreasing $G$-martingales with $K_{0}^{j}=0$ and $K_{T}^{j}\in L_{G}^{\alpha ^{\ast }}(\Omega_{T})$. Then the processus defined by
\begin{equation*}
\int_{0}^{t}X_{s}^{+}dK_{s}^{1}+\int_{0}^{t}X_{s}^{-}dK_{s}^{2}
\end{equation*}%
is also a decreasing $G$-martingale.
\end{Lem}
\section{Large deviations for $G$-SDEs}\label{Sec3}
In this section, we present the large deviations for $G$-SDEs obtained by \citet{Gao2010}. The authors use discrete time approximation to obtain their results.

First, we recall the following notations on large deviations under a sublinear expectation.

Let $(\chi, d)$ be a Polish space. Let $\left(U^{\varepsilon},\;\varepsilon>0\right)$ be a family of measurable maps from $\Omega$ into $(\chi, d)$ and let $\delta(\varepsilon)$, $\varepsilon>0$ be a positive function satisfying $\delta(\varepsilon)\rightarrow 0$ as $\varepsilon\rightarrow 0$.

A nonnegative function $\mathcal{I}$ on $\chi$ is called to be (good) \emph{rate function} if $\{x:\;\mathcal{I}(x)\leq\alpha\}$ (its level set) is (compact) closed for all $0\leq\alpha<\infty$.

$\left\{\widehat{C}(U^{\varepsilon}\in\cdot)\right\}_{\varepsilon>0}$ is said to satisfy large deviation principle with speed $\delta(\varepsilon)$ and with rate function $\mathcal{I}$ if for any measurable closed subset $\mathcal{F}\subset\chi$,
\begin{equation*}\label{UBLD}
    \limsup_{\varepsilon\rightarrow 0}\delta(\varepsilon)\log \widehat{C}\left(U^{\varepsilon}\in \mathcal{F}\right)\leq-\inf_{x\in \mathcal{F}}\mathcal{I}(x),
\end{equation*}
and for any measurable open subset $\mathcal{O}\subset\chi$,
\begin{equation*}\label{LBLD}
    \liminf_{\varepsilon\rightarrow 0}\delta(\varepsilon)\log \widehat{C}\left(U^{\varepsilon}\in \mathcal{O}\right)\geq-\inf_{x\in \mathcal{O}}\mathcal{I}(x).
\end{equation*}
In \citet{Gao2010}, for any $\varepsilon>0$, the authors considered the following random perturbation SDEs driven by $d$-dimensional $G$-Brownian motion $B$
\begin{equation*}
    X_{t}^{x,\varepsilon} =x+\int_{0}^{t}b^{\varepsilon}(X_{r}^{x,\varepsilon})dr+\varepsilon\int_{0}^{t}h^{\varepsilon}(X_{r}^{x,\varepsilon})d\langle B, B\rangle_{r/\varepsilon}+\varepsilon\int_{0}^{t}\sigma^{\varepsilon}(X_{r}^{x,\varepsilon})dB_{r/\varepsilon}
\end{equation*}
where $\langle B, B\rangle$ is treated as a $d\times d$-dimensional vector,
\begin{equation*}
    b^{\varepsilon}=(b_{1}^{\varepsilon}, \ldots, b_{n}^{\varepsilon})^{\tau}:\;\mathbb{R}^{n}\longrightarrow\mathbb{R}^{n},\;
    \sigma^{\varepsilon}=(\sigma_{i, j}^{\varepsilon}):\;\mathbb{R}^{n}\longrightarrow\mathbb{R}^{n\times d}
\end{equation*}
and $h^{\varepsilon}:\;\mathbb{R}^{n}\longrightarrow\mathbb{R}^{n\times d^{2}}$.

Consider the following conditions:
\begin{description}
  \item[(H1)] $b^{\varepsilon}$, $\sigma^{\varepsilon}$ and $h^{\varepsilon}$ are uniformly bounded;
  \item[(H2)] $b^{\varepsilon}$, $\sigma^{\varepsilon}$ and $h^{\varepsilon}$ are uniformly Lipschitz continuous;
  \item[(H3)] $b^{\varepsilon}$, $\sigma^{\varepsilon}$ and $h^{\varepsilon}$ converge uniformly to $b:=b^{0}$, $\sigma:=\sigma^{0}$ and $h:=h^{0}$ respectively.
\end{description}
Let $\mathcal{C}([0, T],\mathbb{R}^{n})$ be the space of $\mathbb{R}^{n}$-valued continuous functions $\varphi$ on $[0, T]$ and $\mathcal{C}_{0}([0, T],\mathbb{R}^{n})$ the space of $\mathbb{R}^{n}$-valued continuous functions $\widetilde{\varphi}$ on $[0, T]$ with $\widetilde{\varphi}_{0}=0$.

Define
\begin{align*}
\mathbb{H}^{d}:=&\Big\{\phi\in\mathcal{C}_{0}([0,T],\mathbb{R}^{d}): \phi\;\text{is absolutely continuous and} \\
&\quad\quad\quad\quad\Vert\phi\Vert_{\mathbb{H}}^{2}:=\int_{0}^{T}\vert\phi^{\prime}(r)\vert^{2}dr<+\infty\Big\}, \\
\mathbb{A}:=&\Big\{\eta=\int_{0}^{t}\eta^{\prime}(r)dr;\; \eta^{\prime}: [0,T]\longrightarrow\mathbb{R}^{d\times d}\; \text{Borel measurable and} \\
&\quad\quad\quad\quad\eta^{\prime}(t)\in\Sigma\;\textrm{ for all } t\in[0, T]\Big\}.
\end{align*}
We recall the following result of a joint large deviation principle for $G$-Brownian motion and its quadratic variation process.
\begin{Thm}\citep[See][p.~2225]{Gao2010}.
$\left\{\widehat{C}\left((\varepsilon B_{t/\varepsilon}, \varepsilon\langle B\rangle_{t/\varepsilon})\mid_{t\in[0, T]}\;\in\cdot\right)\right\}_{\varepsilon>0}$ satisfies large deviation principle with speed $\varepsilon$ and with rate function
\begin{equation*}
J(\phi, \eta)=
\begin{cases}
\frac{1}{2}\int_{0}^{T}\langle \phi^{\prime}(r), (\eta^{\prime}(r))^{-1}\phi^{\prime}(r)\rangle dr,  &  \text{if } (\phi, \eta)\in\mathbb{H}^{d}\times\mathbb{A}, \\
+\infty,  &  \text{otherwise}.
\end{cases}
\end{equation*}
\end{Thm}
For any $(\phi, \eta)\in\mathbb{H}^{d}\times\mathbb{A}$, let $\Psi(\phi,\eta)\in\mathcal{C}([0, T],\mathbb{R}^{n})$ be the unique solution of the following ordinary differential equation (ODE in short)
\begin{eqnarray*}
  \Psi(\phi,\eta)(t) &=& x+\int_{0}^{t}b(\Psi(\phi,\eta)(r))dr+\int_{0}^{t}\sigma(\Psi(\phi,\eta)(r))\phi^{\prime}(r)dr \\
  && +\int_{0}^{t}h(\Psi(\phi,\eta)(r))\eta^{\prime}(r)dr.
\end{eqnarray*}
For $0\leq\alpha<1$ given and $n\geq1$, for each $\psi\in\mathcal{C}_{0}([0, T],\mathbb{R}^{n})$, set
\begin{equation*}
    \Vert\psi\Vert_{\alpha}:=\sup_{s,t\in[0, T]}\frac{\vert\psi(s)-\psi(t)\vert}{\vert s-t\vert^{\alpha}}
\end{equation*}
and
\begin{equation*}
    \mathcal{C}^{\alpha}_{0}([0, T],\mathbb{R}^{n}):=\Big\{\psi\in\mathcal{C}_{0}([0,T],\mathbb{R}^{n}): \lim_{\delta\rightarrow0}\sup_{\vert s-t\vert<\delta}\frac{\vert\psi(s)-\psi(t)\vert}{\vert s-t\vert^{\alpha}}=0,\Vert\psi\Vert_{\alpha}<\infty\Big\}.
\end{equation*}
\begin{Thm}\label{l311}\citep[See][p.~2227]{Gao2010}. Let $0\leq\alpha<1/2$ and let $(H1)$, $(H2)$ and $(H3)$ hold.
Then for any closed subset $\mathcal{F}$ and any open subset $\mathcal{O}$ in $\left(\mathcal{C}^{\alpha}_{0}([0, T],\mathbb{R}^{n}), \Vert\cdot\Vert_{\alpha}\right)$,
\begin{equation*}
    \limsup_{\varepsilon\rightarrow 0}\varepsilon\log\widehat{C}\left((X_{t}^{x, \varepsilon}-x)\mid_{t\in[0, T]}\;\in \mathcal{F}\right)\leq-\inf_{\psi\in \mathcal{F}}I(\psi),
\end{equation*}
and
\begin{equation*}
    \liminf_{\varepsilon\rightarrow 0}\varepsilon\log\widehat{C}\left((X_{t}^{x, \varepsilon}-x)\mid_{t\in[0, T]}\;\in \mathcal{O}\right)\geq-\inf_{\psi\in \mathcal{O}}I(\psi),
\end{equation*}
where
\begin{equation*}\label{}
    I(\psi)=\inf\Big\{J(\phi,\eta):\psi=\Psi(\phi,\eta)-x\Big\}.
\end{equation*}
\end{Thm}
We immediately have the following result.
\begin{Coro}\label{l31} Let $(H1)$, $(H2)$ and $(H3)$ hold. Then for any closed subset $\mathcal{F}$ and any open subset $\mathcal{O}$ in $\mathcal{C}_{0}([0, T],\mathbb{R}^{n})$,
\begin{equation*}
    \limsup_{\varepsilon\rightarrow 0}\varepsilon\log\widehat{C}\left((X_{t}^{x, \varepsilon}-x)\mid_{t\in[0, T]}\;\in \mathcal{F}\right)\leq-\inf_{\widetilde{\varphi}\in \mathcal{F}}\Lambda(\widetilde{\varphi}),
\end{equation*}
and
\begin{equation*}
    \liminf_{\varepsilon\rightarrow 0}\varepsilon\log\widehat{C}\left((X_{t}^{x, \varepsilon}-x)\mid_{t\in[0, T]}\;\in \mathcal{O}\right)\geq-\inf_{\widetilde{\varphi}\in \mathcal{O}}\Lambda(\widetilde{\varphi}),
\end{equation*}
where
\begin{equation*}\label{}
    \Lambda(\widetilde{\varphi})=\inf\Big\{J(\phi,\eta):x+\widetilde{\varphi}=
    \Psi(\phi,\eta)\Big\}.
\end{equation*}
\end{Coro}
\section{Large deviations for $G$-BSDEs with subdifferential operator}\label{Sec4}
We consider the $G$-expectation space $(\Omega_{T}, L_{G}^{1}(\Omega_{T}),\widehat{\mathbb{E}})$ with $\Omega_{T}=\mathcal{C}_{0}([0, T],\mathbb{R})$ and $\overline{\sigma}^{2}=\widehat{\mathbb{E}}(B_{1}^{2})\geq-\widehat{\mathbb{E}}(-B_{1}^{2})=\underline{\sigma}^{2}>0$.
\subsection{Assumptions and problem formulation}
\citet{Yang2017} obtained the existence and uniqueness of the solution of the following backward stochastic differential equation driven by $G$-Brownian motion with subdifferential operator
\begin{equation}\label{GMBSDE0}
\begin{cases}
&-dY_{t}+\partial\Pi(Y_{t})dt\ni f(t,Y_{t},Z_{t})dt-Z_{t}dB_{t}+g(t,Y_{t},Z_{t})d\langle B\rangle_{t}-dK_{t} \\
&Y_{T}=\xi
\end{cases}
\end{equation}
where
\begin{description}
  \item[(A1)] $\Pi$: $\mathbb{R}\longrightarrow(-\infty,\;+\infty]$ is a proper lower semicontinuous (l.s.c. in short) convex function such that $\Pi(y)\geq\Pi(0)=0$, for all $y\in\mathbb{R}$.\\
      Denote
      \begin{eqnarray*}
        \textrm{Dom}(\Pi) &=& \{y\in\mathbb{R}:\;\Pi(y)<\infty\}, \\
        \partial\Pi(y) &=& \{u\in\mathbb{R}:\;\langle u, v-y\rangle+\Pi(y)\leq\Pi(v),\;\forall v\in\mathbb{R}\}, \\
        \textrm{Dom}(\partial\Pi) &=& \{y\in\mathbb{R}:\;\partial\Pi(y)\neq\emptyset\}, \\
        (y,u)\in\textrm{Gr}(\partial\Pi) &\Longleftrightarrow& y\in\textrm{Dom}(\partial\Pi),\;u\in\partial\Pi(y).
      \end{eqnarray*}
      Note that the subdifferential operator $\partial\Pi$: $\mathbb{R}\longrightarrow2^{\mathbb{R}}$ is a maximal monotone operator, that is
      \begin{equation*}\label{}
        \langle y-y', u-u'\rangle\geq0,\;\forall (y,u),(y',u')\in\textrm{Gr}(\partial\Pi).
      \end{equation*}
  \item[(A2)] For any $y,z$, $f(\omega,.,y,z),g(\omega,.,y,z)\in M_{G}^{2}(0,T)$.
  \item[(A3)] The functions $f$: $[0, T]\times\mathbb{R}\times\mathbb{R}\longrightarrow\mathbb{R}$ and $g$: $[0, T]\times\mathbb{R}\times\mathbb{R}\longrightarrow\mathbb{R}$ are continuous and there exists a constant $L>0$ such that for all $t\in[0, T]$, $y,y',z,z'\in\mathbb{R}$,
      \begin{equation*}
        \vert f(t, y, z)-f(t, y', z')\vert + \vert g(t, x, y, z)-g(t, y', z')\vert\leq L(\vert y-y'\vert+\vert z-z'\vert).
      \end{equation*}
\end{description}
\begin{Def}
Let $\xi\in L_{G}^{2}(\Omega_{T})$, the solution of the $G$-MBSDE~\eqref{GMBSDE0} is a quadruple of processes $(Y, Z, K, U)$ such that
\begin{enumerate}
\item $Y\in\mathcal{S}^{2}_{G}(0, T)$, $Z\in H^{2}_{G}(0, T)$, $K$ is a decreasing $G$-martingale with $K_{0}=0$, $K_{T}\in L^{2}_{G}(\Omega_{T})$ and $U\in H^{2}_{G}(0, T)$;
\item
\begin{equation*}
\widehat{\mathbb{E}}\Big(\int_{0}^{T}\Pi(Y_r)dr\Big)<+\infty\,;
\end{equation*}
\item For every $0\leq t\leq T$,
\begin{eqnarray*}
Y_{t} + \int_{t}^{T}U_{r}dr&=&\xi + \int_{t}^{T}f(r,Y_{r},Z_{r})dr + \int_{t}^{T}g(r,Y_{r},Z_{r})d\langle B\rangle_{r} \\
&&-\int_{t}^{T}Z_{r}dB_{r} - (K_{T}-K_{t}),\, q.s.\,;
\end{eqnarray*}
\item $(Y_t, U_t)\in\textrm{Gr}(\partial\Pi)$, q.s. on $\Omega_T\times[0, T]$.
\end{enumerate}
\end{Def}
To establish large deviation principle, we consider the following forward-backward stochastic differential
equation driven by $G$-Brownian motion with subdifferential operator: for every $s\leq t\leq T$, $x\in\mathbb{R}$,
\begin{equation*}\label{}
X_{t}^{s,x,\varepsilon}=x+\int_s^t b(X_{r}^{s,x,\varepsilon})dr+\int_s^t\varepsilon h(X_{r}^{s,x,\varepsilon})d\langle B\rangle_{r}+\int_s^t\varepsilon\sigma(X_{r}^{s,x,\varepsilon})dB_{r}
\end{equation*}
\begin{equation}\label{GMBSDE}
\begin{cases}
&-dY_{t}^{s,x,\varepsilon}+\partial\Pi(Y_{t}^{s,x,\varepsilon})dt\ni f(t,X_{t}^{s,x,\varepsilon},Y_{t}^{s,x,\varepsilon},Z_{t}^{s,x,\varepsilon})dt-Z_{t}^{s,x,\varepsilon}dB_{t} \\
&\qquad\qquad\qquad\qquad\qquad\qquad+g(t,X_{t}^{s,x,\varepsilon},Y_{t}^{s,x,\varepsilon},Z_{t}^{s,x,\varepsilon})d\langle B\rangle_{t}-dK_{t}^{s,x,\varepsilon} \\
&Y_{T}^{s,x,\varepsilon}=\Phi(X_{T}^{s,x,\varepsilon})
\end{cases}
\end{equation}
where $\Pi$ is a proper l.s.c. convex function such that $\Pi(y)\geq\Pi(0)=0$, for all $y\in\mathbb{R}$ and
\begin{equation*}
    b, h, \sigma:\mathbb{R}\longrightarrow\mathbb{R};\;
    \Phi:\mathbb{R}\longrightarrow\mathbb{R};\;f, g:[0, T]\times\mathbb{R}\times\mathbb{R}\times\mathbb{R}\longrightarrow\mathbb{R}
\end{equation*}
are deterministic functions and satisfy the following assumptions:
\begin{description}
  \item[(B1)] $b$, $\sigma$ and $h$ are bounded, i.e., there exists a constant $L>0$ such that
  \begin{equation*}
    \sup_{x\in\mathbb{R}}\max\Big\{\vert b(x)\vert, \vert\sigma(x)\vert, \vert h(x)\vert\Big\}\leq L,
  \end{equation*}
  \item[(B2)] $f$ and $g$ are continuous in $t$;
  \item[(B3)]There exist a constant $L>0$ such that
  \begin{eqnarray*}
     && \vert b(x)-b(x')\vert + \vert h(x)-h(x')\vert + \vert \sigma(x)-\sigma(x')\vert\leq L\vert x-x'\vert, \\
     && \vert \Phi(x)-\Phi(x')\vert\leq L\vert x-x'\vert, \\
     && \vert f(t, x, y, z)-f(t, x', y', z')\vert\leq L(\vert x-x'\vert+\vert y-y'\vert+\vert z-z'\vert), \\
     && \vert g(t, x, y, z)-g(t, x', y', z')\vert\leq L(\vert x-x'\vert+\vert y-y'\vert+\vert z-z'\vert).
  \end{eqnarray*}
\end{description}
It follows from \citet{Yang2017} that, under the assumptions $\bf{(B1)-(B3)}$, the $G$-MBSDE~\eqref{GMBSDE} has a unique solution $\{(Y_{t}^{s,x,\varepsilon},Z_{t}^{s,x,\varepsilon},K_{t}^{s,x,\varepsilon},U_{t}^{s,x,\varepsilon}):
s\leq t\leq T\}$ such that
\begin{itemize}
\item $Y^{s,x,\varepsilon}\in\mathcal{S}^{2}_{G}(0, T)$, $Z^{s,x,\varepsilon}\in H^{2}_{G}(0, T)$, $K^{s,x,\varepsilon}$ is a decreasing $G$-martingale with $K_{s}^{s,x,\varepsilon}=0$, $K_{T}^{s,x,\varepsilon}\in L^{2}_{G}(\Omega_{T})$ and $U^{s,x,\varepsilon}\in H^{2}_{G}(0, T)$;
\item
\begin{equation*}
\widehat{\mathbb{E}}\Big(\int_{s}^{T}\Pi(Y_{r}^{s,x,\varepsilon})dr\Big)<+\infty\,;
\end{equation*}
\item For every $s\leq t\leq T$,
\begin{eqnarray*}
Y_{t}^{s,x,\varepsilon}+\int_{t}^{T}U_{r}^{s,x,\varepsilon}dr&=&\Phi(X_{T}^{s,x,\varepsilon}) + \int_{t}^{T}f(r,X_{r}^{s,x,\varepsilon},Y_{r}^{s,x,\varepsilon},Z_{r}^{s,x,\varepsilon})dr \\
&&+ \int_{t}^{T}g(r,X_{r}^{s,x,\varepsilon},Y_{r}^{s,x,\varepsilon},Z_{r}^{s,x,\varepsilon})d\langle B\rangle_{r} \\
&&-\int_{t}^{T}Z_{r}^{s,x,\varepsilon}dB_{r}-(K_{T}^{s,x,\varepsilon}-K_{t}^{s,x,\varepsilon}),\, q.s.\,;
\end{eqnarray*}
\item $(Y_{t}^{s,x,\varepsilon}, U_{t}^{s,x,\varepsilon})\in\textrm{Gr}(\partial\Pi)$, q.s. on $\Omega_T\times[s, T]$.
\end{itemize}
Our purpose is to study the asymptotic behavior of the family $(Y^{s,x,\varepsilon})_{\varepsilon > 0}$ as $\varepsilon$ goes to $0$.
\subsection{Convergence and large deviation principle for the solution of the backward equation}
We consider the following decoupled forward-backward stochastic differential equation driven by $G$-Brownian motion with subdifferential operator: for every $s\leq t\leq T$,
\begin{equation}\label{syst11}
\begin{cases}
&dX_{t}^{s,x,\varepsilon}=b(X_{t}^{s,x,\varepsilon})dt+\varepsilon h(X_{t}^{s,x,\varepsilon})d\langle B\rangle_{t}+\varepsilon\sigma(X_{t}^{s,x,\varepsilon})dB_{t} \\
&X_{s}^{s,x,\varepsilon}=x
\end{cases}
\end{equation}
\begin{equation}\label{syst12}
\begin{cases}
&-dY_{t}^{s,x,\varepsilon}+U_{t}^{s,x,\varepsilon}dt=f(t,X_{t}^{s,x,\varepsilon},Y_{t}^{s,x,\varepsilon},Z_{t}^{s,x,\varepsilon})dt
-Z_{t}^{s,x,\varepsilon}dB_{t} \\
&\qquad\qquad\qquad\qquad\qquad+g(t,X_{t}^{s,x,\varepsilon},Y_{t}^{s,x,\varepsilon},Z_{t}^{s,x,\varepsilon})d\langle B\rangle_{t}-dK_{t}^{s,x,\varepsilon} \\
&Y_{T}^{s,x,\varepsilon}=\Phi(X_{T}^{s,x,\varepsilon}) \\
&(Y_{t}^{s,x,\varepsilon}, U_{t}^{s,x,\varepsilon})\in\textrm{Gr}(\partial\Pi),\;\widehat{\mathbb{E}}\Big(\int_{s}^{T}\Pi(Y_{r}^{s,x,\varepsilon})dr\Big)<+\infty
\end{cases}
\end{equation}
We also consider the following deterministic system: for every $s\leq t\leq T$,
\begin{equation}\label{syst21}
\begin{cases}
&d\varphi_{t}^{s,x}=b(\varphi_{t}^{s,x})dt \\
&\varphi_{s}^{s,x}=x
\end{cases}
\end{equation}
\begin{equation}\label{syst22}
\begin{cases}
&-d\psi_{t}^{s,x}+U_{t}^{s,x}dt=f(t,\varphi_{t}^{s,x},\psi_{t}^{s,x},0)dt+2G(g(t,\varphi_{t}^{s,x},\psi_{t}^{s,x},0))dt \\
&\psi_{T}^{s,x}=\Phi(\varphi_{T}^{s,x}) \\
&(\psi_{t}^{s,x}, U_{t}^{s,x})\in\textrm{Gr}(\partial\Pi),\;\widehat{\mathbb{E}}\Big(\int_{s}^{T}\Pi(\psi_{r}^{s,x})dr\Big)<+\infty
\end{cases}
\end{equation}
Consider the following decreasing $G$-martingale:
\begin{equation*}
    M_{t}^{s,x}=\int_{s}^{t}g(r, \varphi_{r}^{s,x}, \psi_{r}^{s,x}, 0)d\langle B\rangle_{r}-2\int_{s}^{t}G\left(g(r, \varphi_{r}^{s,x}, \psi_{r}^{s,x}, 0)\right)dr.
\end{equation*}
We have the following result for the solution of the forward equation~\eqref{syst11}.
\begin{Lem}\label{l5}
Let $\textbf{(B1)}$ and $\textbf{(B3)}$ hold. Then
\begin{enumerate}
  \item Let $p\geq 2$. For any $\varepsilon\in(0,1]$, there exists a constant $C_{p}>0$, independent of $\varepsilon$, such that
\begin{equation}\label{eq10}
\widehat{\mathbb{E}}\Big(\sup_{s\leq t\leq T}\vert X_{t}^{s,x,\varepsilon}-\varphi_{t}^{s,x}\vert^{p}\Big)\leq C_{p}\varepsilon^{p}.
\end{equation}
  \item Moreover, $\left\{\widehat{C}\left((X_{t}^{s, x, \varepsilon}-x)\mid_{t\in[s, T]}\;\in\cdot\right)\right\}_{\varepsilon>0}$ satisfies a large deviation principle with speed $\varepsilon$ and with rate function
\begin{equation*}\label{}
    \Lambda(\widetilde{\varphi})=\inf\Big\{J(\phi,\eta):x+\widetilde{\varphi}=
    \widehat{\Psi}(\phi,\eta)\Big\},
\end{equation*}
where $\widehat{\Psi}(\phi,\eta)\in\mathcal{C}([s, T],\mathbb{R})$ be the unique solution of the following ODE
\begin{equation*}
    \widehat{\Psi}(\phi,\eta)(t) = x+\int_{s}^{t}b(\widehat{\Psi}(\phi,\eta)(r))dr.
\end{equation*}
\end{enumerate}
\end{Lem}
We recall a very important result in large deviation theory, used to transfer a LDP from one space to another.
\begin{Lem}\label{l41}(Contraction principle). Let $\{\mu_{\varepsilon}\}_{\varepsilon>0}$ be a family of probability measures that satisfies the large deviation principle with a good rate function $\Lambda$ on a Hausdorff topological space $\chi$, and for $\varepsilon\in(0,1]$, let $f_{\varepsilon}:\;\chi\longrightarrow\Upsilon$ be continuous functions, with $(\Upsilon, d)$ a metric space. Assume that there exists a measurable map $f:\;\chi\longrightarrow\Upsilon$ such that for any compact set $\mathcal{K}\subset\chi$,
\begin{equation*}\label{}
    \limsup_{\varepsilon\rightarrow 0}\sup_{x\in\mathcal{K}}d\left(f_{\varepsilon}(x),\;f(x)\right)=0.
\end{equation*}
Suppose further that $\{\mu_{\varepsilon}\}_{\varepsilon>0}$ is exponentially tight.
Then the family of probability measures $\{\mu_{\varepsilon}\circ f_{\varepsilon}^{-1}\}_{\varepsilon>0}$ satisfies the LDP in $\Upsilon$ with the good rate function
\begin{equation*}
    \Lambda'(y)=\inf\Big\{\Lambda(x): x\in\chi, y=f(x)\Big\}.
\end{equation*}
\end{Lem}
The proofs of Lemmas~\ref{l5} and \ref{l41} can be found in \citet{Dakaou2021}.
\begin{Thm}\label{thm1}
Let $\bf{(B1)-(B3)}$ hold. For any $\varepsilon\in(0,1]$, there exists a constant $C>0$, independent of $\varepsilon$, such that
\begin{equation*}
    \widehat{\mathbb{E}}\Big(\sup_{s\leq t\leq T}\vert Y_{t}^{s,x,\varepsilon}-\psi_{t}^{s,x}\vert^{2}\Big)\leq C\varepsilon^{2}.
\end{equation*}
\end{Thm}
\begin{pf}
We consider the following $G$-BSDE: for every $s\leq t\leq T$, $x\in\mathbb{R}$,
\begin{eqnarray}\label{eq7}
Y_{t}&=&\Phi(\varphi_{T}^{s,x})+\int_{t}^{T}f(r,\varphi_{r}^{s,x},Y_{r},Z_{r})dr-\int_{t}^{T}U_{r}dr
\nonumber \\
&&+\int_{t}^{T}g(r,\varphi_{r}^{s,x},Y_{r},Z_{r})d\langle B\rangle_{r}-\int_{t}^{T}Z_{r}dB_{r}-(K_{T}-K_{t}).
\end{eqnarray}
Thanks to equation \eqref{syst22} and the uniqueness of the solution of the $G$-MBSDEs, it is easy to check that $\{(\psi_{t}^{s,x}, 0, M_{t}^{s,x}, U_{t}^{s,x}): s\leq t\leq T\}$ is the solution of the $G$-MBSDE \eqref{eq7}.
So, we have
\begin{eqnarray*}
Y_{t}^{s,x,\varepsilon }-\psi _{t}^{s,x} &=&\Phi \left(
X_{T}^{s,x,\varepsilon }\right) -\Phi \left( \varphi _{T}^{s,x}\right)
-\int_{t}^{T}\left\{ U_{r}^{s,x,\varepsilon }-U_{r}^{s,x}\right\} dr \\
&&+\int_{t}^{T}\left\{ f\left( r,X_{r}^{s,x,\varepsilon
},Y_{r}^{s,x,\varepsilon },Z_{r}^{s,x,\varepsilon }\right) -f\left(
r,\varphi _{r}^{s,x},\psi _{r}^{s,x},0\right) \right\} dr \\
&&+\int_{t}^{T}\left\{ g\left( r,X_{r}^{s,x,\varepsilon
},Y_{r}^{s,x,\varepsilon },Z_{r}^{s,x,\varepsilon }\right) -g\left(
r,\varphi _{r}^{s,x},\psi _{r}^{s,x},0\right) \right\} d\left\langle
B\right\rangle _{r} \\
&&-\int_{t}^{T}Z_{r}^{s,x,\varepsilon }dB_{r}-\left( K_{T}^{s,x,\varepsilon
}-K_{t}^{s,x,\varepsilon }\right) +\left( M_{T}^{s,x}-M_{t}^{s,x}\right).
\end{eqnarray*}
For $\gamma >0$, by It\^{o}'s formula applied to $e^{\gamma t}\left\vert
Y_{t}^{s,x,\varepsilon }-\psi _{t}^{s,x}\right\vert ^{2}$, we have
\begin{eqnarray*}
&&e^{\gamma t}\left\vert Y_{t}^{s,x,\varepsilon }-\psi _{t}^{s,x}\right\vert
^{2}+\gamma \int_{t}^{T}e^{\gamma r}\left\vert
Y_{r}^{s,x,\varepsilon }-\psi _{r}^{s,x}\right\vert
^{2}dr+\int_{t}^{T}e^{\gamma r}\left\vert Z_{r}^{s,x,\varepsilon
}\right\vert ^{2}d\left\langle B\right\rangle _{r} \\
&=&e^{\gamma T}\left\vert \Phi \left( X_{T}^{s,x,\varepsilon }\right) -\Phi \left(
\varphi _{T}^{s,x}\right) \right\vert^{2}-2\int_{t}^{T}e^{\gamma
r}\left\langle Y_{r}^{s,x,\varepsilon}-\psi _{r}^{s,x}, U_{r}^{s,x,\varepsilon }-U_{r}^{s,x}\right\rangle dr \\
&&+2\int_{t}^{T}e^{\gamma r}\left\langle Y_{r}^{s,x,\varepsilon }-\psi _{r}^{s,x}, f\left( r,X_{r}^{s,x,\varepsilon
},Y_{r}^{s,x,\varepsilon },Z_{r}^{s,x,\varepsilon }\right) -f\left(
r,\varphi _{r}^{s,x},\psi _{r}^{s,x},0\right)\right\rangle dr \\
&&+2\int_{t}^{T}e^{\gamma r}\left\langle Y_{r}^{s,x,\varepsilon }-\psi _{r}^{s,x}, g\left( r,X_{r}^{s,x,\varepsilon
},Y_{r}^{s,x,\varepsilon },Z_{r}^{s,x,\varepsilon }\right) -g\left(
r,\varphi _{r}^{s,x},\psi _{r}^{s,x},0\right)\right\rangle d\left\langle
B\right\rangle_{r} \\
&&-\int_{t}^{T}e^{\gamma r}\left\langle Y_{r}^{s,x,\varepsilon }-\psi _{r}^{s,x}, Z_{r}^{s,x,\varepsilon }\right\rangle
dB_{r}-2\int_{t}^{T}e^{\gamma r}\left\{ Y_{r}^{s,x,\varepsilon }-\psi
_{r}^{s,x}\right\} dK_{r}^{s,x,\varepsilon } \\
&&+2\int_{t}^{T}e^{\gamma r}\left\{ Y_{r}^{s,x,\varepsilon }-\psi _{r}^{s,x}\right\} dM_{r}^{s,x}.
\end{eqnarray*}
Since
\begin{equation*}
\left\langle Y_{r}^{s,x,\varepsilon}-\psi _{r}^{s,x}, U_{r}^{s,x,\varepsilon}-U_{r}^{s,x}\right\rangle\geq 0,
\end{equation*}
and
\begin{equation*}
-2\int_{t}^{T}e^{\gamma r}\left\{ Y_{r}^{s,x,\varepsilon }-\psi
_{r}^{s,x}\right\} ^{+}dK_{r}^{s,x,\varepsilon }-2\int_{t}^{T}e^{\gamma
r}\left\{ Y_{r}^{s,x,\varepsilon }-\psi _{r}^{s,x}\right\}
^{-}dM_{r}^{s,x}\geq 0,
\end{equation*}
by using Young's inequality and Lipschitz conditions of $f$ and $g$, we get
\begin{eqnarray*}
&&e^{\gamma t}\left\vert Y_{t}^{s,x,\varepsilon }-\psi _{t}^{s,x}\right\vert
^{2}+\gamma \int_{t}^{T}e^{\gamma r}\left\vert
Y_{r}^{s,x,\varepsilon }-\psi _{r}^{s,x}\right\vert ^{2}dr+\underline{\sigma
}^{2}\int_{t}^{T}e^{\gamma r}\left\vert Z_{r}^{s,x,\varepsilon }\right\vert
^{2}dr+J_{T}-J_{t} \\
&\leq &e^{\gamma T}\left\vert \Phi \left( X_{T}^{s,x,\varepsilon }\right) -\Phi \left(
\varphi _{T}^{s,x}\right) \right\vert^{2} \\
&&+2\int_{t}^{T}e^{\gamma r}\left\vert
Y_{r}^{s,x,\varepsilon }-\psi _{r}^{s,x}\right\vert\left\vert f\left( r,X_{r}^{s,x,\varepsilon
},Y_{r}^{s,x,\varepsilon },Z_{r}^{s,x,\varepsilon }\right) -f\left(
r,\varphi _{r}^{s,x},\psi _{r}^{s,x},0\right) \right\vert dr \\
&&+2\int_{t}^{T}e^{\gamma r}\left\vert
Y_{r}^{s,x,\varepsilon }-\psi _{r}^{s,x}\right\vert\left\vert g\left( r,X_{r}^{s,x,\varepsilon
},Y_{r}^{s,x,\varepsilon },Z_{r}^{s,x,\varepsilon }\right) -g\left(
r,\varphi _{r}^{s,x},\psi _{r}^{s,x},0\right) \right\vert d\left\langle
B\right\rangle _{r} \\
&\leq &e^{\gamma T}\left\vert \Phi \left( X_{T}^{s,x,\varepsilon }\right) -\Phi \left(
\varphi _{T}^{s,x}\right) \right\vert^{2} \\
&&+2L\left( 1+\overline{\sigma }^{2}\right) \int_{t}^{T}e^{\gamma r}\left\vert
Y_{r}^{s,x,\varepsilon }-\psi _{r}^{s,x}\right\vert\left\{
\left\vert X_{r}^{s,x,\varepsilon }-\varphi _{r}^{s,x}\right\vert
+\left\vert Y_{r}^{s,x,\varepsilon }-\psi _{r}^{s,x}\right\vert +\left\vert
Z_{r}^{s,x,\varepsilon }\right\vert \right\} dr \\
&\leq &e^{\gamma T}\left\vert \Phi \left( X_{T}^{s,x,\varepsilon }\right) -\Phi \left(
\varphi _{T}^{s,x}\right) \right\vert^{2} +L\left( 1+\overline{%
\sigma }^{2}\right) \int_{t}^{T}e^{\gamma r}\left\vert
X_{r}^{s,x,\varepsilon }-\varphi _{r}^{s,x}\right\vert ^{2}dr \\
&&+L\left( 1+\overline{\sigma }^{2}\right) \left(2+\frac{L\left(1+\overline{
\sigma }^{2}\right) }{\underline{\sigma }^{2}}\right) \int_{t}^{T}e^{\gamma
r}\left\vert Y_{r}^{s,x,\varepsilon }-\psi _{r}^{s,x}\right\vert ^{2}dr+%
\underline{\sigma }^{2}\int_{t}^{T}e^{\gamma r}\left\vert
Z_{r}^{s,x,\varepsilon }\right\vert ^{2}dr.
\end{eqnarray*}
where
\begin{eqnarray*}
J_{t} &=&\int_{s}^{t}e^{\gamma r}Z_{r}^{s,x,\varepsilon }\left\{
Y_{r}^{s,x,\varepsilon }-\psi _{r}^{s,x}\right\} dB_{r} \\
&&+2\int_{s}^{t}e^{\gamma r}\left\{ Y_{r}^{s,x,\varepsilon }-\psi
_{r}^{s,x}\right\} ^{+}dK_{r}^{s,x,\varepsilon }+2\int_{s}^{t}e^{\gamma
r}\left\{ Y_{r}^{s,x,\varepsilon }-\psi _{r}^{s,x}\right\} ^{-}dM_{r}^{s,x}.
\end{eqnarray*}
We have, by setting $\gamma =L\left(1+\overline{\sigma }^{2}\right)\left(2+\frac{L\left( 1+\overline{\sigma }^{2}\right) }{\underline{\sigma }^{2}}\right)$,
\begin{eqnarray*}
\left\vert Y_{t}^{s,x,\varepsilon }-\psi _{t}^{s,x}\right\vert
^{2}+J_{T}-J_{t} &\leq &e^{\gamma
T}\left\vert \Phi \left( X_{T}^{s,x,\varepsilon
}\right) -\Phi \left( \varphi _{T}^{s,x}\right) \right\vert^{2} \\
&&\quad\quad+L\left( 1+\overline{\sigma }^{2}\right) \int_{t}^{T}e^{\gamma
r}\left\vert X_{r}^{s,x,\varepsilon }-\varphi _{r}^{s,x}\right\vert ^{2}dr \\
&\leq &e^{\gamma T}\left(L^{2}+L\left( 1+\overline{\sigma }^{2}\right) T\right)
\sup_{s\leq r\leq T}\left\vert X_{r}^{s,x,\varepsilon }-\varphi
_{r}^{s,x}\right\vert ^{2} \\
&\leq & C\sup_{s\leq r\leq T}\left\vert X_{r}^{s,x,\varepsilon }-\varphi
_{r}^{s,x}\right\vert ^{2}.
\end{eqnarray*}
Since $J$ is a $G$-martingale, taking conditional $G$-expectation, we get
\begin{equation*}
\left\vert Y_{t}^{s,x,\varepsilon }-\psi _{t}^{s,x}\right\vert ^{2}\leq C
\widehat{\mathbb{E}}_{t}\left[ \sup_{s\leq r\leq T}\left\vert X_{r}^{s,x,\varepsilon
}-\varphi _{r}^{s,x}\right\vert ^{2}\right].
\end{equation*}
Thus we obtain
\begin{equation*}
\widehat{\mathbb{E}}\left[ \sup_{s\leq t\leq T}\left\vert Y_{t}^{s,x,\varepsilon
}-\psi _{t}^{s,x}\right\vert ^{2}\right] \leq C\widehat{\mathbb{E}}\left[ \sup_{s\leq
r\leq T}\left\vert X_{r}^{s,x,\varepsilon }-\varphi _{r}^{s,x}\right\vert
^{2}\right].
\end{equation*}
So, by virtue of \eqref{eq10}, the proof is complete.\qed
\end{pf}
We have an immediate consequence of Theorem~\ref{thm1}.
\begin{Coro}\label{coro1}
For any $\varepsilon\in(0,1]$ and all $x$ in a compact subset of $\mathbb{R}$, there exists a constant $C>0$, independent of $s$, $x$ and $\varepsilon$, such that
\begin{equation*}\label{eqx1}
\widehat{\mathbb{E}}\Big(\sup_{s\leq t\leq T}\vert Y_{t}^{s,x,\varepsilon}-\psi_{t}^{s,x}\vert^{2}\Big)\leq C\varepsilon^{2}.
\end{equation*}
\end{Coro}
\begin{Lem}
For any $\varepsilon\in(0,1]$, there exists a constant $C>0$, independent of $\varepsilon$, such that
\begin{equation}\label{eqx4}
    \widehat{\mathbb{E}}\Big[\vert K_{T}^{s,x,\varepsilon}\vert^{2}\Big]\leq C.
\end{equation}
\end{Lem}
\begin{Thm}\label{thm3}
Let $\bf{(B1)-(B3)}$ hold. For any $\varepsilon\in(0,1]$, there exists a constant $C>0$, independent of $\varepsilon$, such that
\begin{equation*}
    \widehat{\mathbb{E}}\Big[\int_{s}^{T}\vert Z_{r}^{s,x,\varepsilon}\vert^{2}dr\Big]\leq C\varepsilon.
\end{equation*}
\end{Thm}
\begin{pf}
Similarly as in the proof of Theorem~\ref{thm1}, for $\gamma >0$, by It\^{o}'s formula applied to $e^{\gamma t}\left\vert Y_{t}^{s,x,\varepsilon }-\psi _{t}^{s,x}\right\vert ^{2}$, we have
\begin{eqnarray*}
&&\gamma \int_{s}^{T}e^{\gamma r}\left\vert
Y_{r}^{s,x,\varepsilon }-\psi _{r}^{s,x}\right\vert ^{2}dr+\underline{\sigma}^{2}\int_{s}^{T}e^{\gamma r}\left\vert Z_{r}^{s,x,\varepsilon }\right\vert
^{2}dr+\int_{s}^{T}e^{\gamma r}\left\langle Y_{r}^{s,x,\varepsilon }-\psi _{r}^{s,x}, Z_{r}^{s,x,\varepsilon }\right\rangle
dB_{r} \\
&\leq&e^{\gamma T}\left\vert \Phi \left( X_{T}^{s,x,\varepsilon }\right) -\Phi \left(
\varphi _{T}^{s,x}\right) \right\vert^{2} \\
&&+2\int_{s}^{T}e^{\gamma r}\left\vert
Y_{r}^{s,x,\varepsilon }-\psi _{r}^{s,x}\right\vert\left\vert f\left( r,X_{r}^{s,x,\varepsilon
},Y_{r}^{s,x,\varepsilon },Z_{r}^{s,x,\varepsilon }\right) -f\left(
r,\varphi _{r}^{s,x},\psi _{r}^{s,x},0\right) \right\vert dr \\
&&+2\int_{s}^{T}e^{\gamma r}\left\vert
Y_{r}^{s,x,\varepsilon }-\psi _{r}^{s,x}\right\vert\left\vert g\left( r,X_{r}^{s,x,\varepsilon
},Y_{r}^{s,x,\varepsilon },Z_{r}^{s,x,\varepsilon }\right) -g\left(
r,\varphi _{r}^{s,x},\psi _{r}^{s,x},0\right) \right\vert d\left\langle
B\right\rangle _{r} \\
&&-2\int_{s}^{T}e^{\gamma r}\left\{ Y_{r}^{s,x,\varepsilon }-\psi
_{r}^{s,x}\right\} dK_{r}^{s,x,\varepsilon }+2\int_{s}^{T}e^{\gamma r}\left\{ Y_{r}^{s,x,\varepsilon }-\psi _{r}^{s,x}\right\} dM_{r}^{s,x}.
\end{eqnarray*}
By Lipschitz conditions of $f$ and $g$, we get
\begin{eqnarray*}
&&\gamma \int_{s}^{T}e^{\gamma r}\left\vert
Y_{r}^{s,x,\varepsilon }-\psi _{r}^{s,x}\right\vert ^{2}dr+\underline{\sigma}^{2}\int_{s}^{T}e^{\gamma r}\left\vert Z_{r}^{s,x,\varepsilon }\right\vert
^{2}dr+\int_{s}^{T}e^{\gamma r}\left\langle Y_{r}^{s,x,\varepsilon }-\psi _{r}^{s,x}, Z_{r}^{s,x,\varepsilon }\right\rangle
dB_{r} \\
&\leq &e^{\gamma T}\left\vert \Phi \left( X_{T}^{s,x,\varepsilon }\right) -\Phi \left(
\varphi _{T}^{s,x}\right) \right\vert^{2} \\
&&+2L\left(1+\overline{\sigma }^{2}\right) \int_{s}^{T}e^{\gamma r}\left\vert
Y_{r}^{s,x,\varepsilon }-\psi _{r}^{s,x}\right\vert\left\{
\left\vert X_{r}^{s,x,\varepsilon }-\varphi _{r}^{s,x}\right\vert
+\left\vert Y_{r}^{s,x,\varepsilon }-\psi _{r}^{s,x}\right\vert +\left\vert
Z_{r}^{s,x,\varepsilon }\right\vert \right\} dr \\
&&-2\int_{s}^{T}e^{\gamma r}\left\{ Y_{r}^{s,x,\varepsilon }-\psi
_{r}^{s,x}\right\} dK_{r}^{s,x,\varepsilon }+2\int_{s}^{T}e^{\gamma r}\left\{ Y_{r}^{s,x,\varepsilon }-\psi _{r}^{s,x}\right\} dM_{r}^{s,x}
\end{eqnarray*}
Using Young's inequality, we obtain
\begin{eqnarray*}
&&\gamma \int_{s}^{T}e^{\gamma r}\left\vert
Y_{r}^{s,x,\varepsilon }-\psi _{r}^{s,x}\right\vert ^{2}dr+\underline{\sigma}^{2}\int_{s}^{T}e^{\gamma r}\left\vert Z_{r}^{s,x,\varepsilon }\right\vert
^{2}dr+\int_{s}^{T}e^{\gamma r}\left\langle Y_{r}^{s,x,\varepsilon }-\psi _{r}^{s,x}, Z_{r}^{s,x,\varepsilon }\right\rangle
dB_{r} \\
&\leq &e^{\gamma T}\left\vert \Phi \left( X_{T}^{s,x,\varepsilon }\right) -\Phi \left(
\varphi _{T}^{s,x}\right) \right\vert^{2} +L\left(1+\overline{%
\sigma }^{2}\right) \int_{s}^{T}e^{\gamma r}\left\vert
X_{r}^{s,x,\varepsilon }-\varphi _{r}^{s,x}\right\vert ^{2}dr \\
&&+L\left(1+\overline{\sigma}^{2}\right) \left(3+\frac{4L\left(1+\overline{
\sigma }^{2}\right) }{\underline{\sigma }^{2}}\right) \int_{s}^{T}e^{\gamma
r}\left\vert Y_{r}^{s,x,\varepsilon }-\psi _{r}^{s,x}\right\vert ^{2}dr+%
\frac{\underline{\sigma }^{2}}{4}\int_{s}^{T}e^{\gamma r}\left\vert
Z_{r}^{s,x,\varepsilon }\right\vert ^{2}dr \\
&&-2\int_{s}^{T}e^{\gamma r}\left\{ Y_{r}^{s,x,\varepsilon }-\psi
_{r}^{s,x}\right\}^{+}dK_{r}^{s,x,\varepsilon }+2\int_{s}^{T}e^{\gamma
r}\left\{ Y_{r}^{s,x,\varepsilon }-\psi _{r}^{s,x}\right\}^{+}dM_{r}^{s,x} \\
&&+2\int_{s}^{T}e^{\gamma r}\left\{Y_{r}^{s,x,\varepsilon }-\psi
_{r}^{s,x}\right\}^{-}dK_{r}^{s,x,\varepsilon }-2\int_{s}^{T}e^{\gamma
r}\left\{ Y_{r}^{s,x,\varepsilon }-\psi _{r}^{s,x}\right\}^{-}dM_{r}^{s,x} \\
&\leq &e^{\gamma T}\left\vert \Phi \left( X_{T}^{s,x,\varepsilon }\right) -\Phi \left(
\varphi _{T}^{s,x}\right) \right\vert^{2} +L\left( 1+\overline{%
\sigma }^{2}\right) \int_{s}^{T}e^{\gamma r}\left\vert
X_{r}^{s,x,\varepsilon }-\varphi _{r}^{s,x}\right\vert ^{2}dr \\
&&+L\left(1+\overline{\sigma }^{2}\right) \left(3+\frac{4L\left(1+\overline{
\sigma }^{2}\right) }{\underline{\sigma }^{2}}\right) \int_{s}^{T}e^{\gamma
r}\left\vert Y_{r}^{s,x,\varepsilon }-\psi _{r}^{s,x}\right\vert ^{2}dr+%
\frac{\underline{\sigma }^{2}}{4}\int_{s}^{T}e^{\gamma r}\left\vert
Z_{r}^{s,x,\varepsilon }\right\vert ^{2}dr \\
&&+2e^{\gamma T}\left[(-K_{T}^{s,x,\varepsilon})+(-M_{T}^{s,x})\right]\sup_{s\leq t\leq T}\vert Y_{t}^{s,x,\varepsilon}-\psi_{t}^{s,x}\vert.
\end{eqnarray*}
Thus
\begin{eqnarray*}
&&\frac{3\underline{\sigma}^{2}}{4}\int_{s}^{T}e^{\gamma r}\left\vert Z_{r}^{s,x,\varepsilon }\right\vert
^{2}dr+\int_{s}^{T}e^{\gamma r}\left\langle Y_{r}^{s,x,\varepsilon }-\psi _{r}^{s,x}, Z_{r}^{s,x,\varepsilon }\right\rangle
dB_{r} \\
&\leq &e^{\gamma T}\left\vert \Phi \left( X_{T}^{s,x,\varepsilon }\right) -\Phi \left(
\varphi _{T}^{s,x}\right) \right\vert^{2} +L\left( 1+\overline{%
\sigma }^{2}\right) \int_{s}^{T}e^{\gamma r}\left\vert
X_{r}^{s,x,\varepsilon }-\varphi _{r}^{s,x}\right\vert ^{2}dr \\
&&+\left\{L(1+\overline{\sigma }^{2})\left(3+\frac{4L(1+\overline{
\sigma }^{2})}{\underline{\sigma }^{2}}\right)-\gamma\right\}\int_{s}^{T}e^{\gamma
r}\left\vert Y_{r}^{s,x,\varepsilon }-\psi _{r}^{s,x}\right\vert ^{2}dr \\
&&+2e^{\gamma T}\left[(-K_{T}^{s,x,\varepsilon})+(-M_{T}^{s,x})\right]\sup_{s\leq t\leq T}\vert Y_{t}^{s,x,\varepsilon}-\psi_{t}^{s,x}\vert.
\end{eqnarray*}
We have, by setting $\gamma=L(1+\overline{\sigma }^{2})\left(3+\frac{4L(1+\overline{
\sigma }^{2})}{\underline{\sigma }^{2}}\right)$ and Lipschitz condition of $\Phi$,
\begin{eqnarray*}
&&\frac{3\underline{\sigma}^{2}}{4}\int_{s}^{T}e^{\gamma r}\left\vert Z_{r}^{s,x,\varepsilon }\right\vert
^{2}dr+\int_{s}^{T}e^{\gamma r}\left\langle Y_{r}^{s,x,\varepsilon }-\psi _{r}^{s,x}, Z_{r}^{s,x,\varepsilon }\right\rangle
dB_{r} \\
&\leq &e^{\gamma T}L^2\vert X_{T}^{s,x,\varepsilon }-\varphi _{T}^{s,x}\vert^{2}+LT\left(1+\overline{%
\sigma}^{2}\right)e^{\gamma T}\sup_{s\leq t\leq T}\vert X_{t}^{s,x,\varepsilon }-\varphi_{t}^{s,x}\vert^{2} \\
&&+2e^{\gamma T}\left(\vert K_{T}^{s,x,\varepsilon}\vert+\vert M_{T}^{s,x}\vert\right)\sup_{s\leq t\leq T}\vert Y_{t}^{s,x,\varepsilon}-\psi_{t}^{s,x}\vert.
\end{eqnarray*}
Then
\begin{eqnarray*}
&&\frac{3\underline{\sigma}^{2}}{4}\int_{s}^{T}e^{\gamma r}\left\vert Z_{r}^{s,x,\varepsilon }\right\vert
^{2}dr+\int_{s}^{T}e^{\gamma r}\left\langle Y_{r}^{s,x,\varepsilon }-\psi _{r}^{s,x}, Z_{r}^{s,x,\varepsilon }\right\rangle
dB_{r} \\
&\leq &e^{\gamma T}L\left\{L+T\left(1+\overline{\sigma}^{2}\right)\right\}\sup_{s\leq t\leq T}\vert X_{t}^{s,x,\varepsilon }-\varphi_{t}^{s,x}\vert^{2} \\
&&+2e^{\gamma T}\left(\vert K_{T}^{s,x,\varepsilon}\vert+\vert M_{T}^{s,x}\vert\right)\sup_{s\leq t\leq T}\vert Y_{t}^{s,x,\varepsilon}-\psi_{t}^{s,x}\vert.
\end{eqnarray*}
Therefore
\begin{eqnarray*}
&&\frac{3\underline{\sigma}^{2}}{4}\widehat{\mathbb{E}}\Big[\int_{s}^{T}e^{\gamma r}\vert Z_{r}^{s,x,\varepsilon }\vert^{2}\Big]dr \\
&\leq &e^{\gamma T}L\left\{L+T\left(1+\overline{\sigma}^{2}\right)\right\}\widehat{\mathbb{E}}\Big[\sup_{s\leq t\leq T}\vert X_{t}^{s,x,\varepsilon }-\varphi_{t}^{s,x}\vert^{2}\Big] \\
&&+2e^{\gamma T}\left[\left(\widehat{\mathbb{E}}(\vert K_{T}^{s,x,\varepsilon}\vert^2)\right)^{1/2}+\left(\widehat{\mathbb{E}}(\vert M_{T}^{s,x}\vert^2)\right)^{1/2}\right]\left(\widehat{\mathbb{E}}\Big[\sup_{s\leq t\leq T}\vert Y_{t}^{s,x,\varepsilon}-\psi_{t}^{s,x}\vert^2\Big]\right)^{1/2} \\
&\leq &e^{\gamma T}L\left\{L+T\left(1+\overline{\sigma}^{2}\right)\right\}\widehat{\mathbb{E}}\Big[\sup_{s\leq t\leq T}\vert X_{t}^{s,x,\varepsilon }-\varphi_{t}^{s,x}\vert^{2}\Big] \\
&&+2Ce^{\gamma T}\left[\left(\widehat{\mathbb{E}}(\vert K_{T}^{s,x,\varepsilon}\vert^2)\right)^{1/2}+\left(\widehat{\mathbb{E}}(\vert M_{T}^{s,x}\vert^2)\right)^{1/2}\right]\left(\widehat{\mathbb{E}}\Big[\sup_{s\leq t\leq T}\vert X_{t}^{s,x,\varepsilon}-\varphi_{t}^{s,x}\vert^2\Big]\right)^{1/2}.
\end{eqnarray*}
So, by virtue of \eqref{eq10} and \eqref{eqx4}, the proof is complete. \qed
\end{pf}
\begin{Rem}
As a consequence of Theorems~\ref{thm1} and \ref{thm3}, we get
\begin{equation*}
    \widehat{\mathbb{E}}\Big[\sup_{s\leq t\leq T}\vert Y_{t}^{s,x,\varepsilon}-\psi_{t}^{s,x}\vert^{2}+\int_{s}^{T}\vert Z_{r}^{s,x,\varepsilon}\vert^{2}dr\Big]\leq C\varepsilon,
\end{equation*}
where $C$ is a positive constant.
\end{Rem}
We now want to prove that the process $Y^{s,x,\varepsilon}$ satisfies a LDP. For that reason, we recall the link between Variational Inequality (VI in short) and $G$-MBSDEs. For all $\varepsilon>0$, we consider the following VI
\begin{equation}\label{VI}
\begin{cases}
&\partial_{t}u^{\varepsilon}+\mathcal{L}^{\varepsilon}\left(D^{2}_{x}u^{\varepsilon}, D_{x}u^{\varepsilon}, u^{\varepsilon}, x, t\right)\in\partial\Pi(u^{\varepsilon}(t,x)), \\
&u^{\varepsilon}(T,x)=\Phi(x),\,x\in\mathbb{R}
\end{cases}
\end{equation}
where
\begin{eqnarray*}
  \mathcal{L}^{\varepsilon}\left(D^{2}_{x}u^{\varepsilon}, D_{x}u^{\varepsilon}, u^{\varepsilon}, x, t\right) &=& G\left(H\left(D^{2}_{x}u^{\varepsilon}, D_{x}u^{\varepsilon}, u^{\varepsilon}, x, t\right)\right)+\langle b(x), D_{x}u^{\varepsilon}\rangle \\
  &&+f\left(t, x, u^{\varepsilon}, \langle\varepsilon\sigma(x), D_{x}u^{\varepsilon}\rangle\right),
\end{eqnarray*}
and
\begin{eqnarray*}
  H\left(D^{2}_{x}u^{\varepsilon}, D_{x}u^{\varepsilon}, u^{\varepsilon}, x, t\right) &=& D^{2}_{x}u^{\varepsilon}\varepsilon^{2}\sigma\sigma^{\tau} +2\langle D_{x}u^{\varepsilon}, \varepsilon h(x)\rangle \\
   &&+2g\left(t, x, u^{\varepsilon}, \langle\varepsilon\sigma(x), D_{x}u^{\varepsilon}\rangle\right)
\end{eqnarray*}

Now consider
\begin{equation}\label{eq1}
    u^{\varepsilon}(t,x)=Y_{t}^{t,x,\varepsilon},\; (t,x)\in[0, T]\times\mathbb{R}.
\end{equation}
\begin{equation}\label{eq2}
    u^{0}(t,x)=\psi_{t}^{t,x},\; (t,x)\in[0, T]\times\mathbb{R}.
\end{equation}
In \citet{Yang2017} it is shown that $u^{\varepsilon}$ is a viscosity solution of VI~\eqref{VI} and we have
\begin{equation}\label{l42}
    Y_{t}^{s,x,\varepsilon}=u^{\varepsilon}(t,X_{t}^{s,x,\varepsilon}),\; \forall t\in [s,T].
\end{equation}

Let $\mathcal{C}_{0, s}([s, T],\mathbb{R})$ be the space of $\mathbb{R}$-valued continuous functions $\widetilde{\varphi}$ on $[s, T]$ with $\widetilde{\varphi}_{s}=0$.

Let $s\in[0, T]$ and $\varepsilon\geq 0$. We define the mapping $F^{\varepsilon}:\;\mathcal{C}_{0,s}([s, T],\mathbb{R})\longrightarrow\mathcal{C}([s, T],\mathbb{R})$ by
\begin{equation}\label{eq3}
    F^{\varepsilon}(\widetilde{\varphi})=[t\longmapsto u^{\varepsilon}(t, x+\widetilde{\varphi}_{t})],\; s\leq t\leq T,\;\widetilde{\varphi}\in\mathcal{C}_{0,s}([s, T],\mathbb{R}),
\end{equation}
where $u^{\varepsilon}$ is given by \eqref{eq1} and $u^{0}$ by \eqref{eq2}.

By virtue of \eqref{eq3} and \eqref{l42}, for any $\varepsilon >0$ and all $x\in\mathbb{R}$, we have $Y^{s,x,\varepsilon}=F^{\varepsilon}\left(X^{s,x,\varepsilon}-x\right)$.

We have the following result of large deviations
\begin{Thm} Let $\bf{(B1)-(B3)}$ hold. Then for any closed subset $\mathcal{F}$ and any open subset $\mathcal{O}$ in
$\mathcal{C}([s, T], \mathbb{R})$,
\begin{equation*}
    \limsup_{\varepsilon\rightarrow 0}\varepsilon\log\widehat{C}\left(Y^{s, x, \varepsilon}\in\mathcal{F}\right)\leq-\inf_{\psi\in\mathcal{F}}\Lambda'(\psi),
\end{equation*}
and
\begin{equation*}
    \liminf_{\varepsilon\rightarrow 0}\varepsilon\log\widehat{C}\left(Y^{s, x, \varepsilon}\in\mathcal{O}\right)\geq-\inf_{\psi\in\mathcal{O}}\Lambda'(\psi),
\end{equation*}
where
\begin{equation*}
    \Lambda'(\psi)=\inf\Big\{\Lambda(\widetilde{\varphi}): \psi_{t}=F^{0}(\widetilde{\varphi})(t)=u^{0}(t,x+\widetilde{\varphi}_{t}),t\in[s, T],\widetilde{\varphi}\in\mathcal{C}_{0,s}([s, T], \mathbb{R})\Big\}.
\end{equation*}
\end{Thm}
\begin{pf}
Since the family $\left\{\widehat{C}\left((X_{t}^{s, x, \varepsilon}-x)\mid_{t\in[s, T]}\;\in\cdot\right)\right\}_{\varepsilon>0}$ is exponentially tight (see Lemma~3.4 p.~2235 in \citet{Gao2010}), by virtue of Lemma~\ref{l41} (contraction principle) and Lemma~\ref{l5}, we just need to prove that $F^{\varepsilon}$, $\varepsilon>0$ are continuous and $\{F^{\varepsilon}\}_{\varepsilon>0}$ converges uniformly to $F^{0}$ on every compact subset of $\mathcal{C}_{0,s}([s, T], \mathbb{R})$, as $\varepsilon\rightarrow 0$. Since $u^{\varepsilon}$ is continuous, it is not hard to prove that $F^{\varepsilon}$ is also continuous. The uniform convergence of $\{F^{\varepsilon}\}_{\varepsilon>0}$ is a consequence of Corollary~\ref{coro1}. \qed
\end{pf}
\bibliography{LD_GBSDEs_subdifferential_references}
\end{document}